\documentclass[reqno,a4paper]{amsart} 
\usepackage{amsmath,amssymb,amsxtra,amscd}

\numberwithin{equation}{section} 
\theoremstyle{plain} 
\newtheorem{thm}{Theorem}[section] 
\newtheorem{cor}[thm]{Corollary} 
\newtheorem{lem}[thm]{Lemma} 
\newtheorem{Def}[thm]{Definition} 
\newtheorem{prop}[thm]{Proposition} 
 
\theoremstyle{remark} 
\newtheorem{rem}[thm]{Remark} 

\begin{document} 
\title{The quantum orbit method for generalized flag manifolds} 
\author{Jasper V. Stokman} 
\address{Jasper V. Stokman, 
KdV Institute for Mathematics, Universiteit van Amsterdam, 
Plantage Muidergracht 24, 1018 TV Amsterdam, The Netherlands.} 
\email{jstokman@science.uva.nl} 
\date{June 5th, 2002\\
\indent 2000 {\it Mathematics Subject Classification}. 20G42.} 
\begin{abstract} 
Generalized flag manifolds endowed with the Bruhat-Poisson
bra\-cket are compact Poisson homogeneous spa\-ces, whose
decompositions in symplectic leaves
coincide with their stratifications  
in Schubert cells. In this note it is proved that 
the irreducible $*$-representations of the corresponding quantized 
flag manifolds are also parametrized by their 
Schubert cells. 
An important step is the determination of suitable
algebraic generators of the quantized flag manifolds.
These algebraic generators can be naturally expressed in terms of
quantum  Pl{\"u}cker coordinates. This note complements 
the paper of the author and Dijkhuizen in Commun. Math. Phys. 
{\bf 203} (1999), pp. 297-324, in which these results were
established for a special subclass of generalized flag manifolds.
\end{abstract} 
\maketitle 

\section{Introduction} 

The orbit principle of Kostant and Kirillov predicts a correspondence 
between the irreducible unitary representations of Lie groups 
and the coadjoint orbits of the underlying Lie algebra. 
As a natural generalization of this principle one expects 
a correspondence between the irreducible $*$-representations 
of quantized Poisson homogeneous spaces and their symplectic leaves. 
The key example is due to Soibel'man \cite{S}, 
who showed that the irreducible 
$*$-representations of the standard quantization $\mathbb{C}_q[U]$ 
of the regular functions on a compact, connected, simply connected, 
simple Lie group 
$U$ are parametrized by the symplectic leaves of $U$ (with 
the underlying 
Poisson structure on $U$ given by the Bruhat-Poisson bracket). In this 
note the correspondence is further investigated for generalized 
flag manifolds, which form a substantial subclass of Poisson 
$U$-homogeneous spaces. 

A generalized flag manifold $G/P$, with $G$ the complexification 
of $U$ and $P\subseteq G$ a parabolic subgroup, can be viewed as 
a real $U$-homogeneous space $U/K$, 
with $K\subseteq U$ isomorphic to a compact 
real form of the Levi factor of $P$. The flag manifold $U/K$ is 
a Poisson $U$-homogeneous space with symplectic foliation 
naturally parametrized by the Schubert cells of $G/P$. 
The quantization $\mathbb{C}_q[U/K]$ of $U/K$ can be realized 
as a subalgebra of $\mathbb{C}_q[U]$, defined in terms of 
invariance properties with respect to a suitable quantum subgroup 
$\mathbb{C}_q[K]\subseteq \mathbb{C}_q[U]$. 

It is proved in \cite{SD} that the equivalence classes of the 
irreducible $*$-representations of $\mathbb{C}_q[U/K]$ 
are naturally parametrized by Schubert cells, provided that 
$\mathbb{C}_q[U/K]$ is algebraically generated by certain explicit
products of quantum Pl{\"u}cker coordinates. 
This algebraic assumption was verified case by case for 
special examples of flag manifolds $U/K$, in particular when 
$U/K$ is a Hermitean symmetric space. 
In this note this algebraic assumption is proved 
for arbitrary flag manifold $U/K$ 
using an analogue of the Stone-Weierstrass Theorem for type I 
$C^*$-algebras. 
At the end of the note we extend these 
results to Poisson $U$-homogeneous spaces $U/K^{0}$, where 
$K^{0}$ is the semisimple part of $K$. 

{\it Acknowledgements:} The author is supported by a fellowship 
from the Royal Netherlands Academy of Arts and Sciences (KNAW). 

\section{Pl{\"u}cker coordinates on quantized flag manifolds} 
In this section we compare several different
quantized function algebras on generalized flag manifolds. 
We start by introducing the necessary notations and definitions.

\subsection{Structure theory.}\label{structure} 
Let $\mathfrak{g}$ be a complex simple Lie algebra, 
$\mathfrak{h}\subset\mathfrak{g}$ a Cartan subalgebra and 
$R\subset \mathfrak{h}^*$ the corresponding root system. Let 
$\Sigma=\{\alpha_1,\ldots,\alpha_r\}$ be
a fixed choice of (ordered) 
simple roots and denote $R^+\subset R$ for the corresponding set 
of positive roots. Denote $\mathfrak{b}_{\pm}$ for the Borel subalgebras 
\[\mathfrak{b}_{\pm}=\mathfrak{h}\oplus\bigoplus_{\alpha\in \pm R^+} 
\mathfrak{g}_\alpha, 
\] 
where $\mathfrak{g}_\alpha$ is the root space corresponding to the 
root $\alpha$. 
We denote $\bigl(\cdot,\cdot\bigr)$ for the inner product 
on $\mathfrak{h}^*$ dual to the Killing form on $\mathfrak{g}$, 
and $W$ for the Weyl group of $R$. We use the short hand 
notation $s_i\in W$ for the reflection corresponding to 
the simple root $\alpha_i$. 

Let $G$ be the connected, simply connected, real analytic Lie group with Lie 
algebra $\mathfrak{g}$. 
Let $B_{\pm}\subset G$ be the Borel subgroups corresponding to the 
Lie algebras 
$\mathfrak{b}_\pm$. The standard parabolic subgroups $P_S$ of $G$ 
containing $B_+$ are naturally 
parametrized by subsets $S\subseteq \Sigma$, or 
equivalently, by parabolic sub root systems $R_S=R\cap 
\hbox{span}\{S\}\subseteq R$. 
The Levi factor of the Lie algebra $\mathfrak{p}_S=\hbox{Lie}(P_S)$ 
is given by 
\[\mathfrak{l}_S=\mathfrak{h}\oplus \bigoplus_{\alpha\in R_S} 
\mathfrak{g}_\alpha. 
\] 

Using the identification of $\mathfrak{h}$ with its dual space 
$\mathfrak{h}^*$ via the Killing form, we define 
$H_\alpha\in \mathfrak{h}$ to be the Cartan element 
associated to the coroot $d_\alpha^{-1}\alpha\in \mathfrak{h}^*$ 
($\alpha\in R$), where $d_\alpha=\bigl(\alpha,\alpha\bigr)/2$. 
The real span of the $H_\alpha$'s ($\alpha\in R$) 
is a real form $\mathfrak{h}_0$ of the 
Cartan subalgebra $\mathfrak{h}$. A compact real form 
$\mathfrak{u}$ of $\mathfrak{g}$ can now be chosen in such a way 
that $\mathfrak{k}_S:= 
\mathfrak{p}_S\cap\mathfrak{u}=\mathfrak{l}_S\cap\mathfrak{u}$ 
is a compact real form of $\mathfrak{l}_S$ for any subset $S\subseteq 
\Sigma$. Explicitly, $\mathfrak{u}$ is defined by 
\[\mathfrak{u}=i\mathfrak{h}_0\oplus \bigoplus_{\alpha\in R^+} 
\mathbb{R}(E_\alpha-E_{-\alpha})\oplus \bigoplus_{\alpha\in R^+} 
\mathbb{R}i\bigl(E_\alpha+E_{-\alpha}\bigr), 
\] 
with $E_\alpha\in \mathfrak{g}_\alpha$ root vectors satisfying 
$\lbrack E_\alpha,E_{-\alpha}\rbrack=H_\alpha$, 
$\kappa(E_\alpha,E_{-\alpha})=d_\alpha^{-1}$ and $\lbrack 
E_\alpha,E_\beta\rbrack=c_{\alpha,\beta}E_{\alpha+\beta}$ with 
$c_{\alpha,\beta}\in\mathbb{R}$ whenever $\alpha+\beta\in R$. 
Let $U\subset G$ be the connected Lie subgroup with Lie algebra
$\mathfrak{u}$. The Lie subgroup $U\subset G$ is closed, compact
and simply connected.

We now shortly recall the Bruhat-Poisson structure on $U$ (see 
e.g. \cite{LW} and references 
therein for basic facts on Poisson-Lie groups). 
Define a real Lie subalgebra $\mathfrak{b}\subset \mathfrak{g}$ by 
\[ 
\mathfrak{b}=\mathfrak{h}_0\oplus \bigoplus_{\alpha\in R^+} 
\mathfrak{g}_\alpha. 
\] 
Then $(\mathfrak{g},\mathfrak{u},\mathfrak{b})$ is 
a Manin triple with respect to the imaginary part of the Killing 
form. Thus $\mathfrak{u}$ inherets the structure of a Lie bialgebra.
The associated Poisson-Lie group structure on $U$
is called the Bruhat-Poisson structure.
Poisson-Lie subgroups of $U$ with respect to the Bruhat-Poisson 
bracket on $U$ can be 
classified as follows. Denote $\mathfrak{k}_S^0=\lbrack 
\mathfrak{k}_S,\mathfrak{k}_S\rbrack$ for the semisimple part of the 
reductive Lie algebra $\mathfrak{k}_S$.

\begin{prop} 
Let $K\subseteq U$ be a connected Lie subgroup with Lie algebra 
$\mathfrak{k}$. Then $K\subseteq U$ is a Poisson-Lie subgroup 
with respect to the Bruhat-Poisson bracket on $U$
if and only if 
$\mathfrak{k}_S^0\subseteq \mathfrak{k}\subseteq \mathfrak{k}_S$ 
for some subset $S\subseteq \Sigma$. 
\end{prop} 
\begin{proof} 
First observe that any real subspace $\mathfrak{k}\subseteq
\mathfrak{u}$ satisfying 
$\mathfrak{k}_S^0\subseteq \mathfrak{k}\subseteq \mathfrak{k}_S$ 
is a Lie subalgebra of $\mathfrak{u}$ since 
$\mathfrak{k}_S=\mathfrak{k}_S^0\oplus Z(\mathfrak{k}_S)$ 
with $Z(\mathfrak{k}_S)$ the center of the Lie algebra 
$\mathfrak{k}_S$. 

A connected Lie subgroup 
$K\subseteq U$ is a Poisson-Lie subgroup if and only if the 
orthocomplement $\mathfrak{k}^{\perp}\subseteq\mathfrak{b}$ 
of its Lie algebra 
$\mathfrak{k}$ in $\mathfrak{b}$ is 
an ideal (here the orthocomplement is taken with
respect to the imaginary part of the Killing form). 
This class of Lie subalgebras $\mathfrak{k}\subseteq \mathfrak{u}$ 
is in one-to-one correspondence with 
real Lie subalgebras $\mathfrak{p}\subseteq \mathfrak{g}$ containing 
$\mathfrak{b}$. The correspondence is given by 
$\mathfrak{k}=\mathfrak{p}\cap\mathfrak{u}$ and 
$\mathfrak{p}=\mathfrak{k}\oplus \mathfrak{b}$. 
The real linear spaces 
\[\mathfrak{p}=\xi\oplus\bigoplus_{\alpha\in R^+\cup R_S} 
\mathfrak{g}_\alpha 
\] 
with $S\subseteq \Sigma$ and with real linear spaces 
$\xi$ satisfying 
\[ 
i\mathfrak{h}_{0,S}\oplus\mathfrak{h}_0\subseteq 
\xi\subseteq\mathfrak{h}, \qquad 
\mathfrak{h}_{0,S}=\mathbb{R}-\hbox{span}\{H_\alpha \, | \, \alpha\in 
R_S\}, 
\] 
are all the possible real Lie subalgebras $\mathfrak{p}\subseteq 
\mathfrak{g}$ containing $\mathfrak{b}$ (compare with
the classification of the complex parabolic Lie 
subalgebras of $\mathfrak{g}$). The proposition follows 
by intersecting these Lie 
subalgebras $\mathfrak{p}$ with the compact real form $\mathfrak{u}$. 
\end{proof} 
For a given subset $S\subseteq\Sigma$ we write $K_S\subseteq U$ 
(respectively $K_S^0\subseteq U$) for the connected Lie subgroup 
with Lie algebra $\mathfrak{k}_S$ (respectively $\mathfrak{k}_S^0$). 
Both $K_S\subseteq U$ and $K_S^0\subseteq U$ are closed Poisson-Lie 
subgroups of $U$. The homogeneous spaces $U/K_S$ and $U/K_S^0$ inheret
a natural Poisson structure from $U$ 
(which will also be called the Bruhat-Poisson structure). 
In this note we  study the quantization of the Poisson 
$U$-homogeneous space $U/K_S$ in detail. In section \ref{sectss} 
we formulate the main results for the 
quantization of the Poisson $U$-homogeneous space $U/K_S^0$. 

Observe that $U/K_S$ is isomorphic to the generalized 
flag manifold $G/P_S$ as a real $U$-homogeneous space, 
since $K_S=P_S\cap U$ and $U$ acts transitively on $G/P_S$. 
The symplectic foliation of $U/K_S\simeq G/P_S$ coincides 
with the Schubert cell decomposition of $G/P_S$. The 
Schubert cells of $G/P_S$ are parametrized by the coset space $W/W_S$, 
where $W_S\subseteq W$ is the parabolic subgroup generated by the 
simple reflections $\{ s_i\, | \, i: \alpha_i\in S\}$. The minimal coset 
representatives 
\begin{equation}\label{minimal}
W^S=\{w\in W \, | \, l(ws_\alpha)>l(w),\,\,\,\, \forall\, \alpha\in S 
\} 
\end{equation}
with $l(w)$ the length of the Weyl group element $w\in W$, form a 
complete set of representatives for $W/W_S$. They satisfy 
\[ l(uv)=l(u)+l(v),\qquad \forall\, u\in W^S,\,\, \forall\,v\in W_S. 
\] 

\subsection{Quantum groups.} 
Fix $0<q<1$ and denote $U_q(\mathfrak{g})$ for the Drinfel'd-Jimbo 
quantized universal enveloping algebra over $\mathbb{C}$. 
The algebraic generators are denoted by $K_i^{\pm 1}$ and $X_i^{\pm}$ 
$(i=1,\ldots,r$), where the $K_i^{\pm 1}$ correspond to Cartan 
elements and the $X_i^+$ (respectively $X_i^-$) correspond to the  
root vectors of $\mathfrak{g}$ for the simple
roots $\alpha_i$ (respectively 
$-\alpha_i$). 
For the explicit commutation relations we refer to 
\cite[(3.2)]{SD}. We denote $U_q(\mathfrak{h})$ for the unital, 
commutative subalgebra of $U_q(\mathfrak{g})$ generated by the 
$K_i^{\pm 1}$ ($i=1,\ldots,r$). Recall that $U_q(\mathfrak{g})$ 
is a Hopf-$*$-algebra, with the usual formulas for the counit 
$\epsilon$, comultiplication $\Delta$ and antipode $S$ 
(see \cite[(3.3)]{SD}). Our present choice of 
$*$-structure on $U_q(\mathfrak{g})$ reflects the choice of 
compact real form $\mathfrak{u}$ of $\mathfrak{g}$; 
on the generators of $U_q(\mathfrak{g})$ it is explicitly defined by 
\[ (K_i^{\pm 1})^*=K_i^{\pm 1},\quad 
(X_i^+)^*=q_i^{-1}X_i^-K_i,\quad (X_i^-)^*=q_iK_i^{-1}X_i^+ 
\] 
with $q_i=q^{d_i}$ and $d_i=d_{\alpha_i}=(\alpha_i,\alpha_i)/2$. 

Let $P$ be the weight lattice of the root system $R$. Let 
$P^+$ be the cone of dominant weights and $\{\varpi_1,\ldots, 
\varpi_r\}$ the fundamental weights with respect to 
the fixed choice of (ordered) simple roots 
$\Sigma=\{\alpha_1,\ldots,\alpha_r\}$. 
Denote $V(\lambda)$ for the irreducible finite dimensional 
$U_q(\mathfrak{g})$-module with highest weight $\lambda\in P^+$. 
The weight decomposition of $V(\lambda)$ is written as 
\begin{equation*} 
\begin{split} 
V(\lambda)&=\bigoplus_{\mu\leq\lambda}V(\lambda)_\mu,\\ 
V(\lambda)_\mu&=\{v\in V(\lambda) \, | \, K_i\cdot v= 
q^{(\mu,\alpha_i)}v,\,\,\forall\,i \}, 
\end{split} 
\end{equation*} 
where $\leq$ is the dominance order on $P$ with respect to the
positive roots $R^+$. 
We fix a highest weight vector $0\not=v_\lambda\in V(\lambda)_{\lambda}$ 
for each $\lambda\in P^+$ once and for all. 

The quantized function algebra $\mathbb{C}_q[U]$ is the 
Hopf-$*$-subalgebra of the 
Hopf-$*$-algebra dual of $U_q(\mathfrak{g})$, spanned by the 
matrix coefficients of the finite dimensional $P$-weighted 
$U_q(\mathfrak{g})$-representations (see e.g. \cite[(3.5)]{SD} for 
the definitions of the Hopf-$*$-algebra maps on $\mathbb{C}_q[U]$). 
The analogue of the Peter-Weyl Theorem is the decomposition 
\begin{equation}\label{PW} 
\mathbb{C}_q[U]=\bigoplus_{\lambda\in P^+}W(\lambda), 
\end{equation} 
with $W(\lambda)$ the span of the matrix coefficients of $V(\lambda)$.
Note that $\mathbb{C}_q[U]$ is an $U_q(\mathfrak{g})$-bimodule 
with respect to the left and right regular actions 
\[ (X\cdot a)(Y)=a(YX),\quad (a\cdot X)(Y)=a(XY), 
\] 
where $a\in \mathbb{C}_q[U]$ and $X,Y\in U_q(\mathfrak{g})$. 
The Peter-Weyl decomposition \eqref{PW} then coincides with the 
decomposition of $\mathbb{C}_q[U]$ in simple $U_q(\mathfrak{g})$-bimodules.

\subsection{Quantized flag manifolds and Pl{\"u}cker coordinates.} 
In the remainder of this note we fix an arbitrary 
subset $S\subseteq \Sigma$ of the simple roots. We identify the fixed 
subset $S=\{\alpha_{i_1},\ldots,\alpha_{i_l}\}$ of $\Sigma$ with the 
corresponding subset $\{i_1,\ldots,i_l\}$ of $\{1,\ldots,r\}$. 
We furthermore fix a dominant weight $\Lambda\in P^+$ which is 
supported on $\Sigma\setminus S$, and regular 
dominant with respect to $\Sigma\setminus S$. In other words, 
$\Lambda$ is of the form 
$\Lambda=\sum_{i=1}^rm_i\varpi_i$ with $m_i=0$ for $i\in S$ 
and $m_i>0$ for $i\in \Sigma\setminus S$. 

Let $\widetilde{V}(\Lambda)$ be the finite dimensional irreducible 
$G$-representation with highest weight $\Lambda$ 
and fix a highest weight 
vector $0\not=\widetilde{v}_\Lambda\in 
\widetilde{V}(\Lambda)_\Lambda$. Let $\mathbb{C}[G]$ be the 
algebra of regular functions on $G$. 
Let $F_\Lambda$ be the 
$G$-orbit of the line 
$\mathbb{C}\widetilde{v}_\Lambda$ in the projective space 
$\mathbb{P}(\widetilde{V}(\Lambda))$. 
Then $F_\Lambda$ is naturally isomorphic to $G/P_S$ as a complex 
projective variety. 
The algebra $\mathcal{F}_\Lambda$ 
of regular functions on the affine cone over $F_\Lambda$ 
is the unital subalgebra of $\mathbb{C}[G]$ generated by the matrix 
coefficients $\widetilde{f}_\Lambda=f(\cdot\,\widetilde{v}_\Lambda)$, 
($f\in \widetilde{V}(\Lambda)^*$). In this note we regard 
$\mathcal{F}_\Lambda$ as the algebra of holomorphic polynomials on 
(the affine cone over) $F_\Lambda$. 

The algebra of antiholomorpic 
polynomials on $F_\Lambda$ is defined as follows. 
Let $0\not=\xi_\Lambda\in 
\bigl(\widetilde{V}(\Lambda)^*\bigr)_{-\Lambda}$ 
be a lowest weight vector of the 
$G$-representation $\widetilde{V}(\Lambda)^*$ 
dual to $\widetilde{V}(\Lambda)$. Then the $G$-orbit of the line 
$\mathbb{C}\xi_\Lambda$ in the projective space 
$\mathbb{P}(\widetilde{V}(\Lambda)^*)$ is naturally isomorphic to 
$G/P_S^{opp}$, where $P_S^{opp}$ is the parabolic subgroup of 
$G$ associated to $S\subseteq \Sigma$ which contains 
the opposite Borel subgroup $B_-$. 
We call the unital subalgebra 
$\overline{\mathcal{F}}_\Lambda$ of $\mathbb{C}[G]$ generated 
by the regular functions 
\[G\ni 
g\mapsto \xi_\Lambda(g^{-1}\cdot v),\qquad (v\in \widetilde{V}(\Lambda)), 
\] 
the algebra of antiholomorphic polynomials on $F_\Lambda$. 

We have the following natural quantum 
analogues of the algebras $\mathcal{F}_\Lambda$
and $\overline{\mathcal{F}}_\Lambda$, 
cf. Soibel'man \cite{S1}. Denote $V(\lambda)^*$ for the linear dual of
$V(\lambda)$ ($\lambda\in P^+$). For $f\in V(\lambda)^*$ 
we denote 
$f_\lambda\in W(\lambda)\subset \mathbb{C}_q[U]$ for the matrix 
coefficient 
\[ 
f_\lambda(X)=f\bigl(X\cdot v_\lambda\bigr),\qquad 
\forall\,X\in U_q(\mathfrak{g}). 
\] 

\begin{Def} 
{\bf (a)} The unital subalgebra $\mathcal{F}_\Lambda^q\subseteq 
\mathbb{C}_q[U]$ generated by $f_\Lambda$ \textup{(}$f\in 
V(\Lambda)^*$\textup{)} 
is the quantum algebra of holomorphic polynomials on $F_\Lambda$. 

{\bf (b)} The unital subalgebra 
$\overline{\mathcal{F}}_\Lambda^q\subseteq \mathbb{C}_q[U]$ generated 
by $f_\Lambda^*$ \textup{(}$f\in V(\Lambda)^*$\textup{)} 
is the quantum algebra of antiholomorphic polynomials 
on $F_\Lambda$. 
\end{Def} 
Let $\mathcal{B}$ be a basis of $V(\Lambda)^*$ which is compatible with 
the weight decomposition of $V(\Lambda)$, in the sense that 
each $b\in \mathcal{B}$ is supported on only one weight space 
$V(\Lambda)_\mu$, $\mu=\mu(b)\in P$. 
Then the matrix coefficients $f_\Lambda$ ($f\in \mathcal{B}$) 
are called {\it holomorphic quantum Pl{\"u}cker coordinates} 
on $F_\Lambda$. They algebraically generate the algebra 
$\mathcal{F}_\Lambda^q$. 
The corresponding images $f_\Lambda^*$ under the $*$-involution are called 
{\it antiholomorphic quantum Pl{\"u}cker coordinates} on 
$F_\Lambda$. They algebraically generate 
the algebra $\overline{\mathcal{F}}_\Lambda^q$. 

We next define two closely related $*$-algebras $\mathbb{C}_q[U/K_S]$
and $\mathbb{A}_\Lambda$. 
Let $U_q(\mathfrak{k}_S)\subseteq U_q(\mathfrak{g})$ be the 
unital Hopf-$*$-subalgebra algebraically generated by 
$K_i^{\pm 1}$ ($i=1,\ldots,r$) 
and by $X_j^{\pm}$ ($j\in S$). 
\begin{Def} 
{\bf (a)} 
The $*$-algebra of quantized regular functions on the
Poisson $U$-homogeneous space $U/K_S$ is defined by 
\[\mathbb{C}_q[U/K_S]=\{a\in \mathbb{C}_q[U] \, | \, 
X\cdot a=\epsilon(X)a,\,\,\forall\,X\in U_q(\mathfrak{k}_S)\}. 
\] 
{\bf (b)} Let $\mathbb{A}_{\Lambda}$ be the unital 
$*$-subalgebra of $\mathbb{C}_q[U]$ 
algebraically generated by the elements $f_\Lambda\cdot g_\Lambda^*\in 
\mathbb{C}_q[U]$ for all $f,g\in V(\Lambda)^*$. 
\end{Def} 
\begin{rem} 
The elements $f_\Lambda\cdot g_\Lambda^*\in 
\mathbb{C}_q[U]$ ($f,g\in V(\Lambda)^*$)  
are explicitly given by 
\[(f_\Lambda\cdot g_\Lambda^*)(X)= 
\sum f_\Lambda(X_{(1)})\overline{g_\Lambda(S(X_{(2)})^*)}= 
\sum f(X_{(1)}\cdot v_{\Lambda})\overline{g(S(X_{(2)})^*\cdot 
v_{\Lambda})} 
\] 
for $X\in U_q(\mathfrak{g})$, 
where $\Delta(X)=\sum X_{(1)}\otimes X_{(2)}$. 
\end{rem} 

The $*$-subalgebra $\mathbb{C}_q[U/K_S]\subset \mathbb{C}_q[U]$ 
can alternatively be defined in terms of invariance 
properties with respect to a natural quantum subgroup 
$\mathbb{C}_q[K_S]\subseteq \mathbb{C}_q[U]$ (see 
\cite[(4.4)]{SD}). 
Observe that 
both $\mathbb{A}_\Lambda$ and $\mathbb{C}_q[U/K_S]$ are invariant 
under the right regular $U_q(\mathfrak{g})$-action. 

We next discuss the connections between the algebras $\mathcal{F}_q$, 
$\overline{\mathcal{F}}_q$, 
$\mathbb{A}_\Lambda$ and $\mathbb{C}_q[U/K_S]$. It 
is instructive to consider the classical setup ($q=1$) first. 
Let $\chi_\Lambda$ (respectively $\chi_\Lambda^*$) be 
the character of $P_S$ (respectively $P_S^{opp}$) 
defined by 
\[ p\cdot\widetilde{v}_\Lambda=\chi_\Lambda(p)\widetilde{v}_\Lambda,\qquad 
p^\prime\cdot\xi_\Lambda=\chi_\Lambda^*(p^\prime)\xi_\Lambda 
\] 
for $p\in P_S$ and $p^\prime\in P_S^{opp}$. The classical analogues 
of the holomorphic (respectively antiholomorphic) 
Pl{\"u}cker coordinates are in the space 
$\Gamma(G\times_{P_S}\mathbb{C}_{\chi_\Lambda})$ (respectively 
$\Gamma(G\times_{P_S^{opp}}\mathbb{C}_{\chi_\Lambda^*})$) 
of sections of the line bundle 
$G\times_{P_S}\mathbb{C}_{\chi_\Lambda}$ (respectively 
$G\times_{P_S^{opp}}\mathbb{C}_{\chi_\Lambda^*}$). 
Observe that 
$\mathbb{C}_{\chi_\Lambda}\otimes \mathbb{C}_{\chi_\Lambda^*}$, 
viewed as character of the Levi factor $L_S=P_S\cap P_S^{opp}$, 
is the trivial $L_S$-module. 
We thus obtain a well defined map 
\[ 
\Gamma(G\times_{P_S}\mathbb{C}_{\chi_\Lambda})\times 
\Gamma(G\times_{P_S^{opp}}\mathbb{C}_{\chi_\Lambda^*})\rightarrow 
\Gamma(G/L_S) 
\]
(multiplication map), where $\Gamma(G/L_S)$ is the space of sections of 
the trivial line bundle over $G/L_S$. Hence the classical analogue 
of the algebra $\mathbb{A}_\Lambda$ is contained in the algebra 
of regular functions on $G/L_S$. 
The corresponding statement 
in the quantum setup is the inclusion 
$\mathbb{A}_\Lambda\subseteq \mathbb{C}_q[U/K_S]$, which follows from
\cite[Lem. 4.4]{SD}. 
In this note we prove the following stronger assertion. 
\begin{thm}\label{algthm} 
$\mathbb{A}_\Lambda=\mathbb{C}_q[U/K_S]$. 
\end{thm} 
In other words, a polynomial expression in the 
holomorphic and antiholomorphic quantum Pl{\"u}cker coordinates on 
$F_\Lambda$ is a quantum regular function on $U/K_S$ if and only if 
it has zero weight with respect to the left regular 
$U_q(\mathfrak{h})$-action. Furthermore,  
any quantum regular function on $U/K_S$
can be written in this form.
\begin{rem}\label{AS} 
In \cite{SD} an (a priori) larger 
$*$-algebra $\mathbb{A}_S\supseteq \mathbb{A}_\Lambda$ was studied, 
called the factorized 
$*$-algebra. By definition, $\mathbb{A}_S$ is the span of the matrix 
coefficients $f_\lambda\cdot g_\lambda^*$ ($f,g\in V(\lambda)^*$) 
for all dominant weights 
$\lambda\in P^+$ supported on $\Sigma\setminus 
S$. It was shown in \cite[Lem. 4.4]{SD} that 
$\mathbb{A}_S\subseteq \mathbb{C}_q[U/K_S]$. 
Thus Theorem \ref{algthm} 
implies $\mathbb{A}_S=\mathbb{C}_q[U/K_S]$, which 
was conjectured in \cite[Conj. 4.6]{SD}. 
For a special class of flag manifolds $U/K_S$ the equality 
$\mathbb{A}_S=\mathbb{C}_q[U/K_S]$ has been proved 
by a detailed analysis of the branching rules for 
certain finite dimensional $\mathfrak{g}$-representations
(see \cite[Thm. 4.10]{SD}). This class consists of flag manifolds
$U/K_S$ for which $\mathfrak{p}_S\subset \mathfrak{g}$ 
is maximally parabolic and $(U,K_S)$ is a so-called Gel'fand pair
(the most important examples are the
compact Hermitean symmetric spaces, see \cite[Prop. 4.7]{SD}). 
\end{rem}

The strategy for the proof of Theorem \ref{algthm} 
is as follows. 
We first show that if suitable $C^*$-algebra completions of 
$\mathbb{A}_\Lambda$ and 
$\mathbb{C}_q[U/K_S]$ are the same, then automatically 
$\mathbb{A}_\Lambda=\mathbb{C}_q[U/K_S]$. 
Their $C^*$-algebra completions are shown to coincide by invoking 
a version of the Stone-Weierstrass Theorem for type I 
$C^*$-algebras. For the application of the Stone-Weierstrass type 
Theorem a detailed analysis of the 
irreducible $*$-representations of the two $C^*$-algebras is 
necessary, for which we can resort to (slight modifications of) 
the results in \cite{SD}.

\subsection{Towards a proof of Theorem \ref{algthm}.}\label{Tproof} 
The normalized Haar functional 
on $\mathbb{C}_q[U]$ is the linear functional 
$h: \mathbb{C}_q[U]\rightarrow \mathbb{C}$ 
satisfying $h(1)=1$ (normalization) and satisfying 
$h\equiv 0$ on $W(\lambda)$ 
when $\lambda\in P^+\setminus\{0\}$. Equivalently, 
$h$ is the unique normalized linear functional on $\mathbb{C}_q[U]$ 
satisfying the biinvariance properties 
\[ (h\otimes \hbox{id})\Delta(a)=h(a)1=(\hbox{id}\otimes 
h)\Delta(a),\qquad \forall\, a\in \mathbb{C}_q[U]. 
\] 
A detailed study of Haar functionals in the algebraic framework 
was undertaken by Dijkhuizen and Koornwinder \cite{DK}. They consider 
a special class of Hopf-$*$-algebras called {\it compact quantum group 
algebras}, a class containing $\mathbb{C}_q[U]$ as main example. 
In particular, application of 
\cite[Thm. 3.7]{DK} to $\mathbb{C}_q[U]$ implies that
\[\langle a,b\rangle_h=h(b^*a),\qquad \forall\,a,b\in \mathbb{C}_q[U] 
\]
defines a pre-Hilbert structure on $\mathbb{C}_q[U]$. 
Let $\|\cdot \|_h$ be the corresponding norm 
on $\mathbb{C}_q[U]$. 

\begin{lem}\label{ABweak} 
Let $A\subseteq B\subseteq \mathbb{C}_q[U]$ be subspaces 
invariant under the right regular $U_q(\mathfrak{g})$-action. If 
$A$ is dense in $B$ with respect to the norm $\|\cdot \|_h$, then $A=B$. 
\end{lem} 
\begin{proof} 
Let $A$ and $B$ be as in the lemma, 
and suppose that $A\subseteq B$ is dense but $A\not=B$. 
The Peter-Weyl decomposition \eqref{PW} 
is the isotypical decomposition of $\mathbb{C}_q[U]$ 
under the right regular $U_q(\mathfrak{g})$-action, hence 
\[A=\bigoplus_{\lambda\in P^+}A(\lambda),
\qquad B=\bigoplus_{\lambda\in P^+} B(\lambda), 
\] 
with $A(\lambda)=A\cap W(\lambda)$ and $B(\lambda)=B\cap W(\lambda)$. 
The subspaces $A(\lambda)$ and $B(\lambda)$ ($\lambda\in P^+$)
are finite 
dimensional. Hence, for some dominant weight $\mu\in P^+$,
there exists a nonzero vector $b\in B(\mu)$ which is orthogonal to $A(\mu)$. 
We now use \cite[Prop. 3.4]{DK}, which asserts that 
the Peter-Weyl decomposition \eqref{PW} is a decomposition in 
orthogonal subspaces. This implies that $b$ is orthogonal to 
$A$, contradicting the assumption that $A\subseteq B$ is dense. 
\end{proof}

Let $A$ be a unital $*$-algebra. Let $\mathcal{L}(H)$ be the
$C^*$-algebra of bounded linear operators on the Hilbert space $H$.
A unit preserving $*$-algebra homomorphism $\pi: A\rightarrow
\mathcal{L}(H)$ is called a $*$-representation of $A$.
We say that a $*$-representation 
$(\pi,H)$ of $A$ is irreducible if $H$ and 
$\{0\}$ are the only $\pi(A)$-invariant closed subspaces of $H$. 
Two $*$-representations $(\pi, H)$ and $(\pi^\prime, H^\prime)$ 
of $A$ are said to be equivalent if there exists an isomorphism 
$U: H\rightarrow H^\prime$ of 
Hilbert spaces intertwining the two actions of $A$.

By e.g. \cite[sect. 4]{DK}, we can now define a $C^*$-seminorm on 
the $*$-algebra $\mathbb{C}_q[U]$ by 
\begin{equation}\label{norminfinity} 
\|a\|_\infty={\underset{\pi}{\hbox{sup}}}\bigl(\|\pi(a)\|\bigr),
\qquad a\in \mathbb{C}_q[U], 
\end{equation} 
with $\pi$ running through a complete set of representatives of 
the irreducible $*$-representations of $\mathbb{C}_q[U]$. 
The supremum in \eqref{norminfinity} may as well be taken over 
all $*$-representations $\pi$ of $\mathbb{C}_q[U]$. 

It turns out that 
$\|\cdot\|_\infty$ is a $C^*$-norm on $\mathbb{C}_q[U]$. 
For this it suffices to construct a faithful $*$-representation of 
$\mathbb{C}_q[U]$. Such a faithful $*$-representation can be realized on 
the Hilbert space completion $\mathcal{H}_h$ of the pre-Hilbert space 
($\mathbb{C}_q[U]$, $\langle \cdot,\cdot\rangle_h$) 
by continuous extension of the regular $\mathbb{C}_q[U]$-representation 
\[\pi_h(a)b=ab,\qquad \forall\,a,b\in \mathbb{C}_q[U], 
\] 
see \cite[Lem. 4.2 \& Thm. 4.4]{DK}. 
\begin{lem}\label{estimate} 
For all $a\in \mathbb{C}_q[U]$, we have $\|a\|_h\leq \|a\|_\infty$. 
\end{lem} 
\begin{proof} 
Observe that $\|1\|_h=1$ since the Haar functional $h$ is normalized, 
and that $\|\pi_h(a)\|\leq \|a\|_\infty$ for all $a\in 
\mathbb{C}_q[U]$. It follows that 
\[\|a\|_h=\|\pi_h(a)1\|_h\leq \|\pi_h(a)\|\leq \|a\|_\infty 
\] 
for all $a\in \mathbb{C}_q[U]$. 
\end{proof} 
The following corollary is a direct consequence of Lemma
\ref{ABweak} and Lemma \ref{estimate}.
\begin{cor} 
Let $A\subseteq B\subseteq \mathbb{C}_q[U]$ be subspaces invariant
under the 
right regular $U_q(\mathfrak{g})$-action. If $A$ is dense in $B$ 
with respect to the $C^*$-norm $\|\cdot \|_\infty$, then $A=B$. 
\end{cor} 
The completion $C_q(U)$ of 
$\mathbb{C}_q[U]$ with respect to the $C^*$-norm $\|\cdot 
\|_\infty$ is a unital $C^*$-algebra. The $C^*$-subalgebras 
$\mathcal{A}_\Lambda\subseteq C_q(U)$ and $C_q(U/K_S)\subseteq C_q(U)$ 
are defined to be the closures of $\mathbb{A}_\Lambda$ and 
$\mathbb{C}_q[U/K_S]$ in $C_q(U)$, respectively.
Since $\mathbb{A}_\Lambda\subseteq \mathbb{C}_q[U/K_S]$,
we have the inclusion $\mathcal{A}_\Lambda\subseteq
C_q(U/K_S)$. The previous corollary implies: 
\begin{cor}\label{top} 
If $\mathcal{A}_\Lambda=C_q(U/K_S)$, then 
$\mathbb{A}_\Lambda=\mathbb{C}_q[U/K_S]$. 
\end{cor} 
By Corollary \ref{top} we can resort to the theory 
of $C^*$-algebras for the proof of Theorem \ref{algthm}. 
We first recall some general definitions and facts on $C^*$-algebras, 
starting with the notion of {\it richness}, cf. \cite[Def. 11.1.1]{D}. 
\begin{Def} 
Let $B$ be a $C^*$-algebra. 
An $C^*$-subalgebra $A\subseteq B$ is called rich in $B$ 
if the following two conditions 
are satisfied: 
\begin{enumerate} 
\item[{\bf --}] If $\pi$ is an irreducible 
$*$-representation of $B$, then 
its restriction $\pi|_A$ is an irreducible $*$-representation of $A$. 

\item[{\bf --}] If $\pi$ and $\pi^\prime$ 
are two inequivalent irreducible 
$*$-representations of $B$, then $\pi|_A$ and $\pi^\prime|_A$ 
are inequivalent. 
\end{enumerate} 
\end{Def} 
If $A\subseteq B$ is an inclusion of $*$-algebras, then we denote
the restriction $\pi|_A: A\rightarrow \mathcal{L}(H)$ 
of a $*$-representation $\pi: B\rightarrow \mathcal{L}(H)$
again by $\pi$ if no confusion is possible. 

We need the following analogue of the Stone-Weierstrass 
Theorem (see \cite[Prop. 11.1.6]{D}). 
\begin{thm}\label{SW} 
Let $B$ be a type I $C^*$-algebra. The only rich $C^*$-subalgebra of 
$B$ is $B$ itself. 
\end{thm} 
We apply this theorem 
to the inclusion $\mathcal{A}_\Lambda\subseteq 
C_q(U/K_S)$ of $C^*$-algebras. 
We start with the following elementary lemma. 
\begin{lem}\label{type1} 
The $C^*$-algebra $C_q(U/K_S)$ is of type I. 
\end{lem} 
\begin{proof} 
By \cite[Rem. 5.5]{S}, $C_q(U)$ is of type I.  
Now use that a $C^*$-subalgebra of a type I $C^*$-algebra is again
of type I, see \cite{D}. 
\end{proof} 
\begin{lem}\label{richness} 
The $C^*$-algebra $\mathcal{A}_\Lambda$ is rich in $C_q(U/K_S)$. 
\end{lem} 
Observe that Theorem \ref{algthm} is a 
direct consequence of Lemma \ref{richness}. Indeed, 
Theorem \ref{SW}, Lemma \ref{type1} and Lemma \ref{richness} 
imply $\mathcal{A}_\Lambda=C_q(U/K_S)$, whence 
$\mathbb{A}_\Lambda=\mathbb{C}_q[U/K_S]$ by Corollary \ref{top}.

The proof of Lemma \ref{richness},
which requires a detailed study of the irreducible 
$*$-re\-pre\-sen\-ta\-tions of both the $C^*$-algebras 
$\mathcal{A}_\Lambda$ and $C_q(U/K_S)$, is given in the next
section. 

\section{Irreducible $*$-representations} 
In this section we discuss the representation theory of the 
$*$-algebras $\mathbb{C}_q[U]$, 
$\mathbb{A}_\Lambda$, $\mathbb{C}_q[U/K_S]$ 
and their completions. 

\subsection{The 
fundamental $*$-representation of $\mathbb{C}_q[\hbox{SU}(2)]$.}
\label{fund}
Let $V(=V(\varpi_1))$ be the vector representation of 
$U_q(\hbox{sl}(2,\mathbb{C}))$. Choose an inner product
$\bigl(\cdot,\cdot\bigr)$ on $V$ 
such that $\bigl( X\cdot v,w\bigr)=\bigl(v,X^*\cdot w\bigr)$ for
all $X\in U_q(\hbox{sl}(2,\mathbb{C}))$ and all $v,w\in V$. 
Choose an orthonormal basis $\{e_+,e_-\}$ of 
$V$ such that $e_+$ (respectively $e_-$) is a highest (respectively
lowest) weight vector of $V$. Then the four matrix coefficients 
\[L_{\epsilon\xi}=\bigl(\cdot\,e_\xi,e_\epsilon\bigr),\qquad 
\epsilon,\xi=\pm 
\] 
are algebraic generators of $\mathbb{C}_q[\hbox{SU}(2)]$. 
The Hopf-$*$-algebra structure of $\mathbb{C}_q[\hbox{SU}(2)]$ 
can be completely 
characterized in terms of these generators 
(see e.g. \cite[sect. 3]{SD}). 

The presentation of $\mathbb{C}_q[\hbox{SU}(2)]$ in terms
of the generators $L_{\epsilon\xi}$ can be used to define
explicit $*$-representations of $\mathbb{C}_q[\hbox{SU}(2)]$.
In particular, if we write $\{e_j\}_{j=0}^\infty$ for 
the standard orthonormal basis of
the Hilbert space $l_2(\mathbb{Z}_+)$, then it follows that
\begin{equation*} 
\begin{split} 
\pi_q(L_{++})e_j&=\sqrt{(1-q^{2j})}\,e_{j-1},\qquad 
\pi_q(L_{+-})e_j=-q^{j+1}e_j,\\ 
\pi_q(L_{-+})e_j&=q^je_j,\quad \pi_q(L_{--})e_j= 
\sqrt{(1-q^{2(j+1)})}\,e_{j+1}
\end{split} 
\end{equation*} 
(with the convention that 
$\pi_q(L_{++})e_0=0$) defines an irreducible 
$*$-representation $\pi_q$ 
of $\mathbb{C}_q[\hbox{SU}(2)]$ on $l_2(\mathbb{Z}_+)$
(see e.g. \cite{S} and references therein). Up to tensoring with
one-dimensional $*$-representations, $\pi_q$ is the only infinite
dimensional irreducible $*$-representation of
$\mathbb{C}_q[\hbox{SU}(2)]$, cf. section \ref{U}.

\subsection{Representations of $\mathbb{C}_q[U]$.}\label{U} 

For $i=1,\ldots,r$ let $\phi_i: U_{q_i}(\hbox{sl}(2,\mathbb{C})) 
\rightarrow U_q(\mathfrak{g})$ be the Hopf-$*$-algebra 
embedding defined on the generators by $\phi_i(K_1^{\pm 1})=K_i^{\pm 1}$ 
and $\phi_i(X_1^{\pm})=X_i^{\pm}$. 
The dual map induces a surjective Hopf-$*$-algebra homomorphism 
$\phi_i^*: \mathbb{C}_q[U]\rightarrow \mathbb{C}_{q_i}[\hbox{SU}(2)]$. 
We write $\pi_i=\pi_{q_i}\circ \phi_i^*$ for the corresponding 
lift of $\pi_{q_i}:
\mathbb{C}_{q_i}[\hbox{SU}(2)]\rightarrow
\mathcal{L}(l_2(\mathbb{Z}_+))$ 
to an irreducible $*$-representation of $\mathbb{C}_q[U]$. 

One-dimensional $*$-re\-pre\-sen\-ta\-tions of $\mathbb{C}_q[U]$
can be explicitly constructed as follows.
Let $T=\mathbb{T}^{\times r}$ be the standard $r$-dimensional compact 
torus, where 
$\mathbb{T}=\{z\in \mathbb{C} \, | \, |z|=1 \}$ is the unit circle 
in the complex plane. Denote $t^\mu=t_1^{m_1}\cdots t_r^{m_r}$ for 
$\mu=\sum_{i=1}^rm_i\varpi_i\in P$ 
and $t=(t_1,\ldots,t_r)\in T$. Then 
an one-dimensional $*$-representation  
$\tau_t: \mathbb{C}_q[U]\rightarrow \mathbb{C}$ 
($t\in T$) can uniquely be defined 
by requiring that $\tau_t(a)=a(1)t^\mu$ 
when the quantum regular function $a\in \mathbb{C}_q[U]$ 
has left regular $U_q(\mathfrak{h})$-weight $\mu\in P$ 
(i.e. $K_i\cdot a=q^{(\mu,\alpha_i)}a$ for all $i$). 

Let $\Delta^{(k)}: \mathbb{C}_q[U]\rightarrow 
\mathbb{C}_q[U]^{\otimes (k+1)}$ for $k\in\mathbb{Z}_+$ be the 
iterated coproduct, defined recursively by 
$\Delta^{(0)}=\hbox{id}$ and $\Delta^{(k+1)}=(\Delta\otimes 
\hbox{id}^{\otimes k})\circ\Delta^{(k)}$. 
\begin{thm}[Soibel'man \cite{S}]\label{class} 
Let $w=s_{i_1}s_{i_2}\cdots s_{i_l}$ be a reduced 
expression for the Weyl group element $w\in W$, 
so that $l=l(w)$ is the length of $w$. Then the $*$-representation 
\[ \pi_w=(\pi_{i_1}\otimes \pi_{i_2}\otimes \cdots \otimes 
\pi_{i_l})\circ \Delta^{(l)}: \mathbb{C}_q[U]\rightarrow 
\mathcal{L}\bigl(l_2(\mathbb{Z}_+)^{\otimes l}\bigr) 
\] 
is independent of the reduced expression 
\textup{(}up to equivalence\textup{)}, and 
\[ \{\pi_{w,t}=(\pi_w\otimes \tau_t)\circ\Delta \, | \, w\in W,\,\, t\in 
T \} 
\] 
is a complete set of mutually inequivalent irreducible 
$*$-representations of $\mathbb{C}_q[U]$. 
\end{thm} 
\begin{rem}\label{fact}
An important role in the proof of Theorem \ref{class} is played by
\cite[Thm. 3.1]{S}, which states that $\mathbb{C}_q[U]$ is
algebraically generated by the quantum regular functions $f_\lambda$
and $g_\lambda^*$ for all $f,g\in V(\lambda)$ and all $\lambda\in P^+$.
The proof of Theorem \cite[Thm.
3.1]{S} uses properties of the branching rules for 
tensor products of finite dimensional irreducible
$\mathfrak{g}$-modules when the highest weights become large, 
see e.g. \cite[Prop. 9.2.2]{J}.
Theorem \ref{algthm} should be viewed as the analogue of 
\cite[Thm. 3.1]{S} for the generalized flag manifold $U/K_S$. 
The analogous
statement for the Poisson homogeneous space $U/K_S^0$ is formulated
in Theorem \ref{ss}{\bf (a)}. 
\end{rem}

\subsection{Representations of $\mathcal{A}_\Lambda$.}\label{repA} 

The $*$-algebra
$\mathbb{A}_\Lambda$ has the convenient property that it is
defined in terms of explicit algebraic generators involving
quantum Pl{\"u}cker coordinates. This property enables one to analyze the
$*$-re\-pre\-sen\-ta\-tions of $\mathbb{A}_\Lambda$ 
in a similar manner as for $\mathbb{C}_q[U]$, cf. Remark \ref{fact}.
This analysis was carried out in \cite[sect. 6]{SD}, 
but then for the a priori larger $*$-algebra 
$\mathbb{A}_S$ (see Remark \ref{AS}). 
Straightforward adjustments of the proof of \cite[Thm 6.13]{SD} 
though show that 
$\{ \bigl(\pi_w, l_2(\mathbb{Z}_+)^{\otimes l(w)}\bigr) 
\,\, | \,\, w\in W^S 
\}$ is a complete set of mutually inequivalent irreducible 
$*$-representations of $\mathbb{A}_\Lambda$, where $W^S\subseteq W$ is given
by \eqref{minimal}. 

The irreducible $*$-representation $\pi_w$ ($w\in W^S$) of
$\mathbb{A}_\Lambda$ is the restriction to 
$\mathbb{A}_\Lambda$ of an $*$-representation of
$\mathbb{C}_q[U]$. Hence it extends by continuity to an 
irreducible $*$-representation of 
$\mathcal{A}_\Lambda$, which we again denote by $\pi_w$.
On the other hand, any irreducible $*$-representation of 
$\mathcal{A}_\Lambda$ restricts to an irreducible $*$-representation
of $\mathbb{A}_\Lambda$. We conclude that: 
\begin{cor}\label{repS1} 
$\{ \pi_w \,\, | \,\, w\in W^S \}$ is 
a complete set of mutually inequivalent irreducible 
$*$-representations of the $C^*$-algebra $\mathcal{A}_\Lambda$. 
\end{cor} 

\subsection{Representations of $\mathbb{C}_q[U/K_S]$.}\label{repUK} 
In this section we first establish the analogue of Corollary 
\ref{repS1} for $C_q(U/K_S)$. Since we (a priori) do not have a 
complete set of algebraic generators of the 
$*$-algebra $\mathbb{C}_q[U/K_S]$, we cannot analyze
its $*$-representations in the same manner  as for
$\mathbb{A}_\Lambda$. The alternative approach is by
analyzing the ireducible $*$-representations of the $C^*$-algebra
$C_q(U/K_S)$ directly. For this, we first observe that the inclusion 
$\mathbb{A}_\Lambda\subseteq \mathbb{C}_q[U/K_S]$ and the results 
in section \ref{repA} imply that 
$\pi_w$ ($w\in W^S$) are mutually inequivalent irreducible 
$*$-representations of $\mathbb{C}_q[U/K_S]$ and of 
$C_q(U/K_S)$. 
Furthermore, by \cite[Prop. 5.7]{SD} we have 
\[\pi_{w,t}(a)=\pi_u(a)\otimes \hbox{id}^{\otimes l(v)},\qquad 
\forall\,a\in C_q(U/K_S) 
\] 
for $t\in T$ and $w=uv\in W$
with $u\in W^S$ and $v\in W_S$. Since every
irreducible $*$-representation of $C_q(U/K_S)$ 
appears as irreducible component of $\pi|_{C_q(U/K_S)}$ for some 
irreducible $*$-representation $\pi$ of $C_q(U)$
(see e.g. \cite[Prop. 2.10.2]{D}), we conclude (cf. \cite[Thm. 5.9]{SD}):
\begin{cor}\label{repS2} 
$\{ \pi_w \,\, | \,\, w\in W^S\}$ 
is a complete set of mutually inequivalent irreducible 
$*$-representations of the $C^*$-algebra $C_q(U/K_S)$. 
\end{cor}
Corollary \ref{repS1} and Corollary \ref{repS2} show that 
$\mathcal{A}_\Lambda$ is rich in $C_q(U/K_S)$. Thus
Lemma \ref{richness} is proved, which in turn implies the validity
of Theorem \ref{algthm} (see section \ref{Tproof}).

In particular, the classification of the irreducible 
$*$-representations of $\mathbb{A}_\Lambda$ (see section \ref{repA}) 
gives the following theorem. 
\begin{thm}\label{repSnew} 
The set 
\[ \{ \bigl(\pi_w, l_2(\mathbb{Z}_+)^{\otimes l(w)}\bigr) 
\,\, | \,\, w\in W^S 
\} 
\] 
is a complete set of mutually inequivalent irreducible 
$*$-representations of the quantized function 
algebra $\mathbb{C}_q[U/K_S]$. 
\end{thm} 
\begin{rem} 
{\bf (a)} 
Theorem \ref{repSnew} is not a direct consequence of 
Corollary \ref{repS2}, since it is a priori not clear that every 
irreducible $*$-representation of $\mathbb{C}_q[U/K_S]$ can be 
continuously extended to a $*$-representation of $C_q(U/K_S)$. 

{\bf (b)} Theorem \ref{repSnew} fits nicely in the philosophy 
of the Kostant-Kirillov orbit method. Indeed, the irreducible 
$*$-representations of $\mathbb{C}_q[U/K_S]$ are parametrized by 
the coset representatives $W/W_S$, which in turn also parametrize 
the symplectic leaves of the underlying Poisson $U$-homogeneous 
space $U/K_S$, see section \ref{structure} and \cite[sect. 2]{SD}. 
\end{rem} 

\section{The Poisson $U$-homogeneous space $U/K_S^{0}$.}\label{sectss} 

In this section we apply the methods of the present note to 
the Poisson $U$-homogeneous space $U/K_S^0$. 

Denote $U_q(\mathfrak{k}_S^{0})$ for the unital Hopf-$*$-subalgebra 
of $U_q(\mathfrak{g})$ generated by $K_i^{\pm 1}, 
X_i^{\pm}$ ($i\in S$). We call  
\[ 
\mathbb{C}_q[U/K_S^{0}]=\{ a\in\mathbb{C}_q[U] \,\, | \,\, 
X\cdot a=\epsilon(X)a,\,\, \forall X\in U_q(\mathfrak{k}_S^{0}) 
\} \subseteq \mathbb{C}_q[U]
\] 
the $*$-algebra of quantized regular functions on the Poisson 
$U$-homogeneous space $U/K_S^{0}$. Note that 
$\mathbb{C}_q[U/K_S^{0}]$ is invariant under the right regular 
$U_q(\mathfrak{g})$-action and under the left regular 
$U_q(\mathfrak{h})$-action. Observe furthermore that 
\[ 
\mathbb{C}_q[U/K_S]=\{a\in \mathbb{C}_q[U/K_S^0]\,\, | \,\, X\cdot 
a=\epsilon(X)a \quad \forall X\in U_q(\mathfrak{h}) \}. 
\] 
\begin{thm}\label{ss} 
{\bf (a)}  
$\mathbb{C}_q[U/K_S^{0}]$ is algebraically generated by 
the matrix coefficients $f_k$ and $g_l^*$ for all $f\in V(\varpi_k)^*$, 
$g\in V(\varpi_l)^*$ and all $k,l\in \Sigma\setminus S$.\\ 
{\bf (b)} A complete set of mu\-tual\-ly i\-neq\-ui\-va\-lent 
ir\-re\-du\-ci\-ble $*$-re\-pre\-sen\-ta\-tions 
of the $*$-algebra $\mathbb{C}_q[U/K_S^{0}]$ is given by 
\[ \{ \bigl(\pi_{w,t}, l_2(\mathbb{Z}_+)^{\otimes l(w)}\bigr) 
\, | \, w\in W^S,\,\, t\in T_{\Sigma\setminus S} \}, 
\] 
with $T_{\Sigma\setminus S}\subseteq T$ the sub-torus 
\[ T_{\Sigma\setminus S}= 
\{ t=(t_1,\ldots,t_r)\in T \, | \, t_i=1\, \hbox{when }\, i\in 
S\}. 
\] 
\end{thm} 
\begin{rem} 
{\bf (a)} The special case $S=\emptyset$ of Theorem \ref{ss} recovers
Soibel'man's \cite{S} results for $\mathbb{C}_q[U]$. 
The irreducible $*$-representations for the special case that
$U/K_S^{0}$ is the Stiefel manifold $\hbox{SU}(n)/\hbox{SU}(l)$ ($l<n$)
were described before by
Podkolzin and Vainerman \cite{PV}.

{\bf (b)} The unital $*$-algebra $\mathbb{A}_\Lambda^{0}$ 
algebraically generated by the holomorphic and antiholomorphic 
quantum Pl{\"u}cker coordinates $f_\Lambda$ and $g_\Lambda^*$ ($f,g\in 
\mathcal{B}$) is in general properly contained in 
$\mathbb{C}_q[U/K_S^{0}]$. This can for instance be verified by 
comparing the weights occurring in $\mathbb{A}_\Lambda^0$ and 
in $\mathbb{C}_q[U/K_S^0]$ 
under the left regular $U_q(\mathfrak{h})$-action.

\end{rem} 
\bibliographystyle{amsplain} 
 
\end{document}